%% file: main.tex
\documentclass[12pt, a4paper]{article}

\include{mycommands_sp}

\author{Simon-Christian Klein\footnote{Institute for Partial Differential Equations, TU Braunschweig, Braunschweig, Germany,
	\href{siklein@tu-bs.de}{siklein@tu-bs.de}
}}

\title{Essentially Non-oscillatory Spectral Volume Methods}

\begin{document}
	\maketitle

	\begin{abstract}
		A new Essentially Non-oscillatory (ENO) recovery algorithm is developed and tested in a Finite Volume method. 
		The construction is hinged on a reformulation of the reconstruction as the solution to a variational problem. The sign property of the classical ENO algorithm is expressed as restrictions on the admissible set of solutions to this variational problem.
		In conjunction with an educated guessing algorithm for possible locations of discontinuities an ENO reconstruction algorithm without divided differences or smoothness indicators is constructed. No tunable parameters exist apart from the desired order and stencil width. The desired order is in principle arbitrary, but growing stencils are needed. 
		While classical ENO methods consider all connected stencils that surround a cell under consideration the proposed recovery method uses a fixed stencil, simplifying efficient high order implementations.
	\end{abstract}
		
	\section{Introduction} \label{sec:intro}
	\input{intro}
	\section{Theory} \label{sec:theory}
	\input{theory}

	\section{Numerical Tests} \label{sec:NT}
	\input{tests}
	\section{Conclusion} \label{sec:concl}
	\input{concl}

	\section{Conflict of Interest Disclosure}
	The author has no relevant financial or non-financial interests to disclose.
	
	\section{Bibliography}
	\bibliographystyle{plainnat}
	\bibliography{lit}
\end{document}

%% file: mycommands_sp.tex
\usepackage{graphicx}
\usepackage[english]{babel}
\usepackage{amsmath}
\usepackage{caption}
\usepackage{comment}
\usepackage{subcaption}
\usepackage{amssymb}
\usepackage{mathtools}
\usepackage{mathptmx} 
\usepackage{xcolor}
\usepackage{tikz}
\usepackage[square,numbers,sort]{natbib}
\usepackage{dsfont}
\usepackage{hyperref}
\usepackage{pgfplots}
\pgfplotsset{width=\textwidth,compat=1.9}
\newcommand{\epswidth}{0.9\textwidth}

\newcommand{\derive}[2] {{\frac{\partial {#1} }{\partial {#2}}}}
\newcommand{\derd}[2] {{\frac{\vd #1 }{\vd #2}}}

\newcommand{\R}[0]{\mathbb{R}}

\newcommand{\abs}[1]{\left | #1\right| }

\newcommand{\norm}[1]{\left \|#1  \right \|}

\newcommand{\set}[2]{ \left \{ #1 ~\middle|~  #2 \right\}}
\newcommand{\sset}[1]{ \left \{ #1 \right \}}

\newcommand{\Leb}{\mathrm{L}}
\DeclareMathOperator{\sign}{\mathrm{sign}}

\newcommand{\vd}{ \mathrm{d}}
\newcommand{\intd}{\,\mathrm{d}}
\newcommand{\e}{\mathrm{e}}
\renewcommand{\epsilon}{\varepsilon}
\renewcommand{\phi}{\varphi}
\newcommand{\skp}[2]{\left \langle {#1}, {#2} \right \rangle}

\newcommand{\txfrac}[2]{#1 \slash #2}

\DeclareMathOperator{\argmin}{argmin}
\DeclareMathOperator{\st}{\mathrm{s.t.}}
\DeclareMathOperator{\bigO}{\mathcal{O}}

\newcommand{\Recov}{\mathcal{R}}
\newcommand{\Eop}{\mathcal{E}}

\newcommand{\Aop}{\mathcal{A}}
\newcommand{\Aset}{\mathtt{A}}
\newcommand{\Nop}{\mathtt{N}}

\newcommand{\TV}{\mathrm{TV}}
\newcommand{\vwave}{{\rotatebox[origin=c]{90}{$\sim$}}}
\newcommand{\minmod}{\mathrm{minmod}}
\newcommand{\smi}[1]{\vwave #1 \vwave}

%% file: intro.tex
Hyperbolic systems of conservation laws
\begin{equation} \label{eq:HCL}
	\sum_{i=1}^n \derive{f_i(u)}{x_i} + \derive{u}{t} = 0
	\end{equation}
for $u(x, t): \R^n \times \R_{\geq 0} \mapsto \R^k$ are defined by $n$ flux functions $f_i: \R^k \mapsto \R^k$. Depending on the number of conserved variables $k$ and the flux functions are those partial differential equation systems able to model flows in subsonic and supersonic regimes \cite{CF1948SFSW}, waves in shallow water \cite{CMMOT2024highorder} and magnetic effects in plasma \cite{WJS2023Provably}. 
In nearly all applications equation \eqref{eq:HCL} is not solved exactly as there exist only a small number of known exact solutions.
Instead, solutions are approximated using numerical techniques, especially using Finite Difference (FD) \cite{GKO2013Time}, Finite Volume (FV) \cite{MS2007FV}, Discontinuous Galerkin (DG) \cite{CJST1997Advanced} and Continuous Galerkin (CG) \cite{CJST1997Advanced} methods.
We will concentrate on Finite Volume methods in one space dimension in this publication. Our domain will be decomposed into several disjoint subdomains with edges $x_{k + \txfrac 1 2}$, called Finite Volumes or cells. Let us denote the cell centers as $x_k = \txfrac{\left(x_{k-\txfrac 1 2} + x_{k + \txfrac 1 2}\right)}{2}$. Integrating our conservation law over one space-time prism $\Omega_k \times [t_1, t_2] = \left[x_{k - \txfrac 1 2}, x_{k + \txfrac 1 2}\right] \times [t_1, t_2]$ we find
\[
	\begin{aligned}
	&\int_{t_1}^{t_2}\int_{\Omega_k} \derive{ u(x, t)}{t} + \derive{f(u(x, t))}{x}\intd x \intd t = 0 \\
	\implies &\int_{\Omega_k} u(x, t_2)  \intd x = \int_{\Omega_k}u(x, t_1) \intd t +  \int_{t_1}^{t_2} f\left(u\left(x_{k - \txfrac 1 2}, t\right)\right) - f\left(u\left(x_{k + \txfrac 1 2}, t\right)\right) \intd t.
	\end{aligned}
\]
Let us define \[
\mu(\Omega_k) =  \int_{\Omega_k} 1 \intd x = x_{k+\txfrac 1 2} - x_{k-\txfrac 1 2} 
\] for the volume of $\Omega_k$ and the abbreviation $u_k(t)$ for the average value of $u$ in cell $k$ at time $t$
\[
	u_k(t) = \frac{1}{\mu(\Omega_k)} \int_{\Omega_k} u(x, t) \intd x.
\]
Using these abreviations follows the evolution equation 
\[
	u_k\left(t_2\right) = u_k\left(t_1\right) + \frac 1 {\mu\left(\Omega_k \right)} \int_{t_1}^{t_2} {f\left(u\left(x_{k - \txfrac 1 2}, t\right)\right) - f\left(u\left(x_{k + \txfrac 1 2}, t\right)\right)} \intd t
\]
for the average values of $u$ on the cells in time. The values of $u$ in the cells can be recovered at $t_1$ by assuming constant functions in each cell. Therefore exist left and right limits of this recovered  $u$ at every cell boundary. The time integral over the flux at $x_{k + \txfrac 1 2}$ with a discontinuous $u$ can be evaluated using an exact Riemann solver, as pioneered by Godunov \cite{Godunov1959A}.
Exact Riemann solvers are not always available and in any case costly. One often resides to approximate Riemann solvers like the HLL family of fluxes \cite{HLL1983On}.
We will now outline how this basic construction was generalized to higher order shock capturing schemes.
Dividing by $t_2 - t_1$ and taking the limit $t_2 \to t_1$ one finds the semidiscrete version
\[
	\derd{u_k}{t} =  \frac {f\left(u\left(x_{k - \txfrac 1 2}, t\right)\right) - f\left(u\left(x_{k + \txfrac 1 2}, t\right)\right)}{\mu\left(\Omega_k \right)},
\]
that can be integrated in time using a (higher order) Runge-Kutta method \cite{SO1988,SO1989}.
This construction would be still only first order accurate, as the recovered function $u$ is only first order exact.
A missing piece for an evaluation of the right hand side with high accuracy is a higher order recovery.
One of the main difficulties in the design of such recovery methods are discontinuities in the solution $u$ called shocks and contact discontinuities, depending on their type \cite{CF1948SFSW,Smoller1994Shock,Toro2009Riemann,BF2003Generalized}. 
Specialized methods were developed to overcome the difficulty that high order polynomial expansions oscillate in the presence of such discontinuities. 

 Van Leer introduced the MUSCL schemes. A piecewise affine linear function is reconstructed in every cell to recover point values of $u$ at cell edges with up to second order accuracy \cite{Leer1973TowardsI,Leer1974TowardsII,Leer1977TowardsIII,Leer1977TowardsIV,Leer1979TowardsV}. The key element is the selection of the recovered slope $s_k$ of the function $u$ in each cell. Given the average values of $3$ cells $u_{k-1}, u_k, u_{k+1}$ two different first order polynomials for cell $k$ can be determined. A left one $p_{k - \txfrac 1 2}(x) = u_k + s_{k - \txfrac 1 2}(x - x_k)$ that has the average value $u_{k-1}$ on the cell $\Omega_{k-1}$, and the average value $u_{k}$ on the cell  $\Omega_{k}$. The second  polynomial $p_{k - \txfrac 1 2}(x) = u_k + s_{k + \txfrac 1 2}(x - x_k)$ can be selected to have the correct average values on the central and right cell, i.e. a stencil that extends to the right.
 Central to Van-Leers approach is the usage of these two candidate slopes $s_{k-\txfrac 1 2}, s_{k + \txfrac 1 2}$ to define the slope $s_k$ of the recovery. The selected slope should keep the resulting function monotone if the average values are monotone, avoid creating new maxima and confine the total variation. Possible choices include the minmod limiter function \cite{Sweby1984High}
 \[
 	\begin{aligned}
 	s_k &= \minmod\left(s_{k - \txfrac 1 2}, s_{k + \txfrac 1 2} \right)\\
 		 	&=  \begin{cases} 0 							& 	s_{k-\txfrac 1 2} s_{k + \txfrac 1 2} < 0 \\
 				 		\begin{cases}
 				 		s_{k - \txfrac 1 2} & \abs{s_{k- \txfrac 1 2}} \leq \abs{s_{k + \txfrac 1 2}}  \\
 						s_{k + \txfrac 1 2} &   \abs{s_{k + \txfrac 1 2}} \leq \abs{s_{k - \txfrac 1 2}} 
 						\end{cases}					& 	\text{else}.
 		\end{cases}
 		\end{aligned}
 \]
 A recovery using the minmod function chooses a zero slope if the slope changes between the left and right polynomial. The accuracy degrades whenever this occurs. Therefore second order is only achieved away from local maxima.

A more advanced construction can be carried out with piecewise parabolic polynomials leading to the Piecewise Parabolic Method (PPM) \cite{CW1984PPM}, achieving up to third order accuracy.

A similar method leads to in principle arbitrary accuracy. Classical Essentialy-Non-Oscillatory (ENO) schems recover a polynomial of degree $q$ that reproduces the correct average values on $q$ cells around a target cell and on the target cell itself \cite{HEOC1987UniformlyIII}. The crucial step lies in the selection of the correct stencil to construct this recovery polynomial. In classical ENO methods the stencil for cell $\Omega_k$ is chosen iteratively. Assume we already found a stencil $S_k^q = \sset{\Omega_k, \Omega_{k + 1}, \ldots, \Omega_{k+q}}$ for a $q$-th order polynomial around $\Omega_k$. Then the classical ENO stencil selection selects the stencil for the polynomial degree $q+1$ by considering two new candidate stencils. Each enlarged by one cell, one to the left and one to the right of $S_k^q$. The highest divided differences of the two candidate polynomials are calculated. 
If the highest divided differences of the left stencil have lower absolute value the left stencil is used and vice-versa.

Numerical experiments indicate that the resulting recovery polynomial $p$ satisfies 
\[
 	\TV(p) \leq \TV(u) + \bigO((\Delta x)^q)
\]
also for non-smooth functions $u$, as indicated in the name \cite{HEOC1987UniformlyIII}. Further, one can show that the \emph{sign property}
\[
	\sign\left( p_{k+1}\left(x_{k + \txfrac 1 2}\right) -  p_{k}\left(x_{k + \txfrac 1 2}\right)\right) = \sign \left (u_{k + 1} - u_k \right)
\]
 holds \cite{FMT2013ENO}. This property allows the construction of entropy dissipative ENO schemes \cite{FMT2012TECNO}.
Further, the jump heights are bounded
\[
	\abs{p_{k+1}\left(x_{k + \txfrac 1 2} \right) - p_k\left(x_{k + \txfrac 1 2}\right)} \leq C_p \abs{u_{k + 1} - u_k}
\]
by the jumps in the average values. This shows that the jumps vanish at extrema as the jumps in average values vanish at extrema, which is essential for the stiffness of the resulting ODE system
\cite{Klein2023EAR}.
ENO methods of degree $q$ can select any connected stencil of size $q+1$ that includes $\Omega_k$ to recover a function in cell $\Omega_k$. Therefore, the average values of $2q+1$ cells are needed to achieve order $q+1$. Selecting a stencil for a high degree in more than one dimension under all connected stencils becomes prohibitive as the number of possible stencils grows extremely fast \cite{Sonar1997ENO}.

Weighted Essentially Non Oscillatory (WENO) methods remove the iterative stencil selection \cite{LOC1994WENO}. Also, they allow a higher order of accuracy for the same number of visited cells.
Smoothness indicators $\smi \cdot \to \R$ are designed to measure if a polynomial is smooth or tries to approximate a discontinuous function.
For a given degree $q$ and cell $k$ a total of $q+1$ recovery polynomials $p^l_k(x)$ are recovered, one for every possible stencil $S_k^l$, $l=1, \ldots, q+1$.
Their smoothness is measured using the smoothness indicator. Coefficients $c^l_k\left (p^1_k, p^2_k, \ldots, p^{q+1}_k \right) \geq 0$ are calculated to define the recovery polynomial via a convex combination
\[
	p_k(x) = c^1_k p^1_k + c^2_k p^2_k + \ldots + c_k^{q+1} p^{r+1}_k = \sum_{l=1}^{q+1} c_k^l p^l_k.
\]
The coefficients $c_k^l$ are functions of the indicators that should tend to zero for non-smooth polynomials, while at the same time converging to a set of base weights $c^l$ for smooth polynomials.
Evaluating the polynomial $p_k$ for a smooth function can result in a higher accuracy than a classical ENO method would achieve if the base weights are chosen correctly.

The classic WENO method does not satisfy a sign property but one can construct modified WENO methods that satisfy the sign property \cite{FR2016SPWENO}.

All presented recoveries are suited for shock-capturing calculations and were originally constructed in one space dimension to reach medium to high orders of accuracy. 
Their recovery operators are non-linear and their nonlinearities are tailored to suppress oscillations.
Spectral Volume (SV) schemes on the other hand can reach in principle arbitrary orders and use a linear recovery \cite{Wang2002SVI,Wang2002SVII,Wang2004SVIII,Wang2004SVIV,Wang2006SVV}.
Assume we subdivided every cell $\Omega_k$ into $q + 1$ subcells \linebreak $\omega_{k, l} = \left[x_{k, l-\txfrac 1 2}, x_{k, l+\txfrac 1 2}\right]$. One can define the averaging operator
\[
	\Aop: \Leb^2(\Omega_k) \to \R^{q+1}, \quad \phi \mapsto \left( \frac{1}{\mu(\omega_l)} \int_{\omega_l} \phi(x) \intd x \right)_l
\]
that maps a function $\phi$ onto the corresponding average values of the function on the subcells of cell $\Omega_k$.
Let $V \subset \Leb^2(\Omega_k)$ be a subspace of $\Leb^2$ with a continuous point evaluation
\[
	\Eop(x): V \to \R,\quad \phi \mapsto \phi(x),
										\]
that maps functions from $V$ to their values at a point $x$ in $\Omega_k$. When $V$ is a $(q+1)$-dimensional subspace, for example the space of all polynomials of degree less than or equal $q$, a $q + 1$ function basis $B = \sset{\phi_1, \phi_2, \phi_3, \ldots, \phi_{q+1}}$ of $V$ exists. The Operator $\Aop$ can be represented on this basis by a $(q+1) \times (q+ 1)$ matrix 
\[
	[\Aop]_B
	= \left  ( \begin{array}{c | c | c | c} 
		 &	& 	& \\
		\Aop \phi_1 & \Aop \phi_2 & \ldots & \Aop \phi_{r+1} \\
		& & & \\
	\end{array} \right )
\]
that contains the images of the basis under the operator $\Aop$. For non-degenerate subcell distributions the restriction of $\Aop$ to this space is invertible, and we can define a linear recovery for the point values via
\[
	\Recov(x) = \Eop(x) \circ \Aop^{-1}.
\]
This recovery operator is the basis for the SV method \cite{Wang2002SVI,Wang2002SVII,Wang2004SVIII,Wang2004SVIV,Wang2006SVV}. Assume for the recovered point values exist left and right sided limits at the subcell edges
\[
	u_{k, l + \txfrac 1 2}^- = \lim_{x \uparrow x_{k, l+\txfrac 1 2}} \Recov(x) u, \quad u_{k, l + \txfrac 1 2}^+ = \lim_{x \downarrow x_{k, l+\txfrac 1 2}} \Recov(x) u. 
\] 
The resulting point recovery can be used in a standard FV method for the average values $u_{k,l}$ of the subcells
\[
\begin{aligned}
	f_{k, l+\frac 1 2} &= f\left(u_{k, l+\txfrac 1 2}^-, u_{k, l+\txfrac 1 2}^+ \right), \quad
	\derd {u_{k,l}}{t} &= \frac{f_{k, l-\txfrac 1 2} -  f_{k, l + \txfrac 1 2}}{\mu \left(\omega_{k, l}\right)}.
	\end{aligned}
\]
This approach is the classical SV method. The cost of integrating the flux over the surfaces of the spectral volumes grows with the dimension and order of the method \cite{ZS2005An}. An additional interpolation step of the flux can reduce the flux integration costs \cite{HWL2008Efficient}, but we will not use this interpolation as it only enlarges the efficiency for multi-dimensional methods.  
Classical SV methods need additional robustness enhancements for shock-capturing calculations as the linear recovery used has no monotonicity, TV, ENO or sign properties. Instead, common limiters from the DG framework are used \cite{ZS2005An,XLS2009Hierachical,LGG2016A}. 
The rest of this publication will design a new \emph{non-linear} recovery operator $\Recov$ that aims to deliver superior performance for shock calculations and removes the need for limiters and similar techniques. This new recovery enforces usual properties of (W)ENO recoveries by restricting the set of admissible recovery functions.

%% file: theory.tex
The next subsections describe our recovery algorithm.
A theoretical formulation is given in subsection \ref{ss:selectiontheory}. Our practical implementation is described in subsection \ref{ss:practicalrecovery}, followed by an explanation how the resulting quadratic variational problems are handled numerically in subsection \ref{ss:quadprog}.

\subsection{Essentially Non-oscillatory Recovery using Subset Selection} \label{ss:selectiontheory}
We propose a new Essentially Non-oscillatory shock capturing scheme. Central is the following nonlinear recovery that takes a set of averages $u_{k, l}$ on a set of cells $\Omega_k$ subdivided into subcells $\omega_{k, l}$. 
From now on the entire recovery is described for the subcells $\omega_{k, l}$ of one single macrocell $\Omega_k$. For clarity the index of the macrocell will be suppressed from now on.
The procedure has the following steps:
\begin{enumerate}
\item Select a (convex) subset $C \subset V$ of target functions for our recovery. This convex subset depends on the average values at hand and introduces a non-linearity.

\item Find $\phi \in C$ minimizing the vector $v \in \R^{r+1}$ defined by
	\[
	v_l = \frac{1}{\mu(\omega_{k,l})}\int_{\omega_l} \phi(x) \intd x - u_{k, l}
	\]
	in the $2$-norm. This can be written using our averaging operator as
	\[
		\phi = \argmin_{\phi \in C} \norm{\Aop \phi - \Aop u}_2.
	\]
	This variational problem uses a general convex set $C$ as restriction. The resulting dependency of $\phi$ on $u$ is in general nonlinear.

\end{enumerate}
\noindent
While this method seems obvious the details are crucial and differentiate our method from existing \guillemetleft{}least squares\guillemetright{}  approaches:
\begin{itemize}
	\item We select $C$ depending on the solution in every cell and every time-step to reflect the nature of the solution. We hope to achieve essentially non-oscillatory behaviour using the selection of a suitable set $C$.
	\item The set $C$ will contain discontinuous functions. 
	\item Classical SV schemes select $C = V$ as polynomials of a specified degree. $V$ is therefore a subspace of $\Leb^2$ and the recovery a linear operator. This implies that classical SV schemes can not be high order accurate and non-oscillatory \cite{Godunov1959A}.
	\item Our sets $C$ will be in general no complete subspaces. We can use this to enforce classical ENO properties like the sign property.
	\item Our recovery operator can not be a precomputed matrix, instead a non-linear solver is needed and we will explain the efficient solution of the resulting non-linear problem.
\end{itemize}
Let us make some observations before we explain our proposed $C$. We will from now on assume that $\dim V \leq q+1$ holds and that the averaging operator $\Aop$ has full rank on the selected subcells and basis functions. This implies that if $C = V$ holds a least squares residual solution is sought and that this solution is unique. In the case $q+1 < \dim V$ the function $\phi$ would be non-unique.
This rules out the simple but ineffective decision to use for example the first $q+1$ monomials \emph{and at the same time} a set of $q$ functions with jumps between the $q+1$ cells as basis functions. In that case our information about the sought function $u$ would be insufficient.

Instead, we will consider the first $k$ monomials and $l$ functions that jump between cells as vectorspace $V$. The total number must satisfy the restriction $k + l \leq q + 1$. 

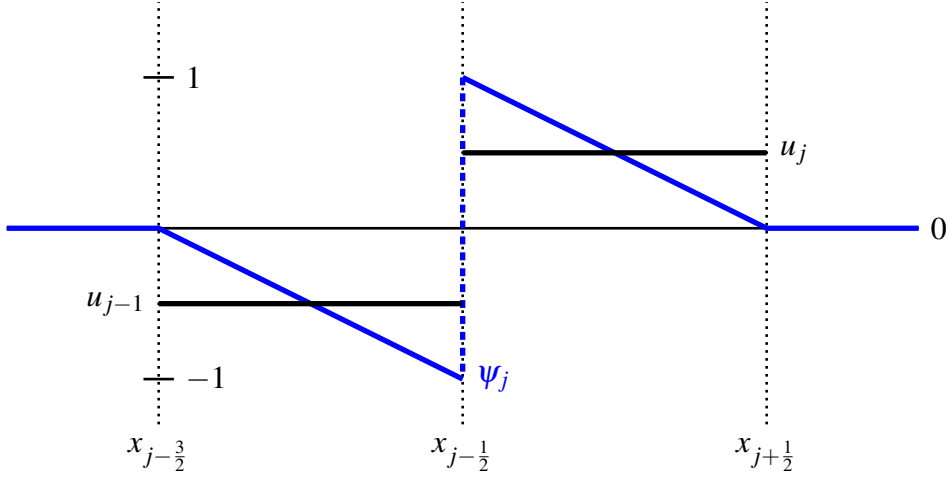
\begin{figure}
	\centering
	\begin{tikzpicture}[scale=2.0]
		\draw [line width=1pt](-3.0, 0.0) -- (3.0, 0.0) node [anchor=west] {$0$};
		\draw (-2.0, -1.5) node {$x_{j - \frac 3 2}$};
		\draw [dotted, line width=1pt] (-2.0, -1.3)  -- (-2.0, 1.5);
		\draw (0.0, -1.5) node {$x_{j - \frac 1 2}$};
		\draw [dotted, line width=1pt] (0.0, -1.3) -- (0.0, 1.5);
		\draw (2.0, -1.5) node {$x_{j + \frac 1 2}$};
		\draw [dotted, line width=1pt] (2.0, -1.3) -- (2.0, 1.5);
		\draw [blue, line width=2pt] (-3.0, 0.0) -- (-2.0, 0.0) -- (0.0, -1.0) node [anchor=west] {$\psi_j$};
		\draw [blue, line width=2pt] (0.0, 1.0) -- (2.0, 0.0) -- (3.0, 0.0);
		\draw [dashed, blue, line width=2pt] (0.0, 1.0) -- (0.0, -1.0);
		\draw [line width=2pt] (-2.0, -0.5) node [anchor = east] {$u_{j-1}$} -- (0.0, -0.5);
		\draw [line width=2pt] (0.0, 0.5) -- (2.0, 0.5) node[anchor = west] {$u_j$};
		\draw [line width=1pt] (-2.1, 1.0) -- (-1.9, 1.0) node [anchor = west] {$1$};
		\draw [line width=1pt] (-2.1, -1.0) -- (-1.9, -1.0) node [anchor = west] {$-1$};
		\end{tikzpicture}
\caption{The prototype inter cell jump function $\Psi_{j-\frac 1 2}(x)$, centerd at the cell edge between cells $j-1$ and $j$. The support of the function is limited to the cells $j-1$ and $j$. The average values of the function in cell $j-1$ is $-\txfrac 1 2$, while it is $\txfrac 1 2$ in cell $j$. The function jumps from $-1$ to $1$ at the edge $x_{j-\frac 1 2}$.}
	\label{fig:ProtoJmp}	
	\end{figure}
The prototypical jump function is
\[
	\psi_{j - \frac 1 2}(x) = \begin{cases} 0 & x \leq x_{j- \frac 3 2}, \\
							-\frac{x - x_{j-\frac 3 2}}{x_{j - \frac 1 2} - x_{j - \frac 3 2}}& x_{j-\frac 3 2} \leq x \leq x_{j- \frac 1 2}, \\
							-\frac{x - x_{j + \frac 1 2}}{x_{j+\frac 1 2} - x_{j - \frac 1 2}}& x_{j-\frac 1 2} \leq x \leq x_{j + \frac 1 2}, \\
							0 & x_{j + \frac 1 2} \leq x,
							\end{cases}
\]
as shown in figure \ref{fig:ProtoJmp}.
This function only has support in the cells $j-1$ and $j$, jumps from a height of $-1$ to $1$ between both and has average $-\txfrac 1 2$ and $\txfrac 1 2$ in the two cells adjacent to the jump. Our complete convex function set $C$ selection algorithm using $k$ continuous and $l$ discontinuous basis functions is given by the following procedure
\begin{enumerate}
	\item Calculate the jumps of the averages 
		\[
			I_{j + \frac 1 2} = u_{j + 1} - u_{j}
			\]
			between all subcells $j = 1, \ldots, q$ that are part of the macrocell.
	\item Sort the jumps by decreasing absolute magnitude $\abs{I_{j_1}} > \abs{I_{j_2}} > \abs{I_{j_3}} > \ldots > \abs{I_{j_{q}}}$.
	\item Select the $l$ edges between two subcells with the highest jumps as the $l$ positions for jumps included as discontinouos basis functions in the basis of $V$.
	\item Add the first $k$ monomials $\phi_i$ to $V$.
	\item Construct the convex subset $C \subset V$ out of the complete vector space $V$ by restriction. Let 
		\begin{equation} \label{eq:lincomb}
			u(x, t) = \sum_{i=1}^k c_i(t) \phi_i(x) + \sum_{i =1}^l d_i(t) \psi_{j_i}(x)
		\end{equation}
		be a linear combination for a function $u \in C$.
			We enforce the sign property by demanding
			\[
			\sign{I_{j_i}} = \sign{d_i}.
			\]
		Our set $C$ is therefore given by
		\[
			C = \set{ \sum_{i=1}^k c_i \phi_i + \sum_{i =1}^l d_i \psi_{j_i}}{c_i \in \R \wedge d_i \in \R \wedge \sign{I_{j_i}} = \sign{d_i} }.
		\]
\end{enumerate}
	Together with a suitable non-linear optimization problem solver this $C$ selection algorithm allows us to recover piecewise smooth functions $u$ from averages.
	These respect the sign property \cite{FMT2013ENO,FR2016SPWENO}. Several other constraints could be placed, like positive density and pressure.

\subsection{Practical Implementation of the Recovery and Cellsize Distributions} \label{ss:practicalrecovery}
	After we defined our recovery in mathematical terms let us now point out how this recovery can be implemented.
	We use orthogonal polynomials instead of the monomials as a basis $\phi_i$ for the continuous part.
	The Legendre polynomials satisfy \cite{HW2008DG,CHQZ2006Spectral}
	\[
		\skp{\phi_i}{\phi_j} =	\int_{-1}^1 \phi_i(x) \phi_j(x) \intd x = \delta_{ij}.
	\]
	Point evaluations yield
	\[
		\sum_k \phi_i(x_k) \phi_j(x_k) \approx \skp{\phi_i}{\phi_j}.
	\]
	This implies that vectors of point evaluations of Legendre polynomials are not as colinear as point evaluations of monomials. This reduces the otherwise exploding condition numbers of Vandermonde matrices
	\[
		V_{ij} = \phi_{j}(x_i), \quad \Aop_{ij} = \frac{1}{x_{i+\frac 1 2} - x_{i - \frac 1 2}} \int_{x_{i-\frac 1 2}}^{x_{i + \frac 1 2}} \phi_j(x) \intd x
	\]
	and average computation \cite{HW2008DG}.
	Given these two operators we can write the linear recovery operator of a classic spectral volume method $\Recov$ as
	\[
		\Recov = V \Aop^{-1}. 
	\]
	This operator maps the averages to a linear combination of basis functions using the least squares inverse of $\Aop$ and evaluates this linear combination at nodal values using the matrix $V$.
	
	It is known from DG and SV methods that non-equidistant nodes or non-equisized cells should be used in the implementation of high order
	schemes \cite{Wang2004SVIII,HW2008DG,Hesthaven1998Electrostatics,CW2015Comparison}.
	Otherwise, the norms of operators extrapolating nodal values to nodal values at other coordinates or average values to nodal values at cell boundaries explode.
	These operator norms are called Lebesgue constants in classical polynomial interpolation \cite[Proposition 4.1]{DL1993Constructive}.
Following \cite{Wang2004SVIII}, we will use the Chebyshev nodes of the second kind 
\[
	x_i = \cos\left(\frac{i-1}{q+1} \pi \right), \quad i=1, \ldots, q+2
\]
as the $q+2$ cell boundaries of our $q+1$ cells. While this cell distribution is not optimal the Lebesgue constant is less than one third higher than for an optimal layout up to order $8$ \cite{Wang2004SVIII}.
After we explained our cell layout we will now explain our solution procedure for the nonlinear part of the recovery - determining the coefficients of the combined basis in equation \eqref{eq:lincomb}.
\subsection{Convex Projections using the Active Set Method} \label{ss:quadprog}
The main stage in our recovery is the projection of our information, i.e. the average values of the sought function $u(x, t)$, onto the convex set of admissible functions $C$. We decided to use a specialized active set method \cite[section 16.4]{NW1999NO} that is combined with the method of conjugate gradients \cite{HS1952CG} for this purpose.
This bounds the computational complexity to be equal to several conjugate gradient linear system solutions per spectral element. 
The worst case complexity for a CG solution of a linear system are $N$ iterations. This amounts to $\bigO(N^3)$ calculations, the same as matrix inversion via QR decomposition or Gaussian elimination. Often, the CG method converges significantly faster, i.e. after significantly less than $N$ iterations.

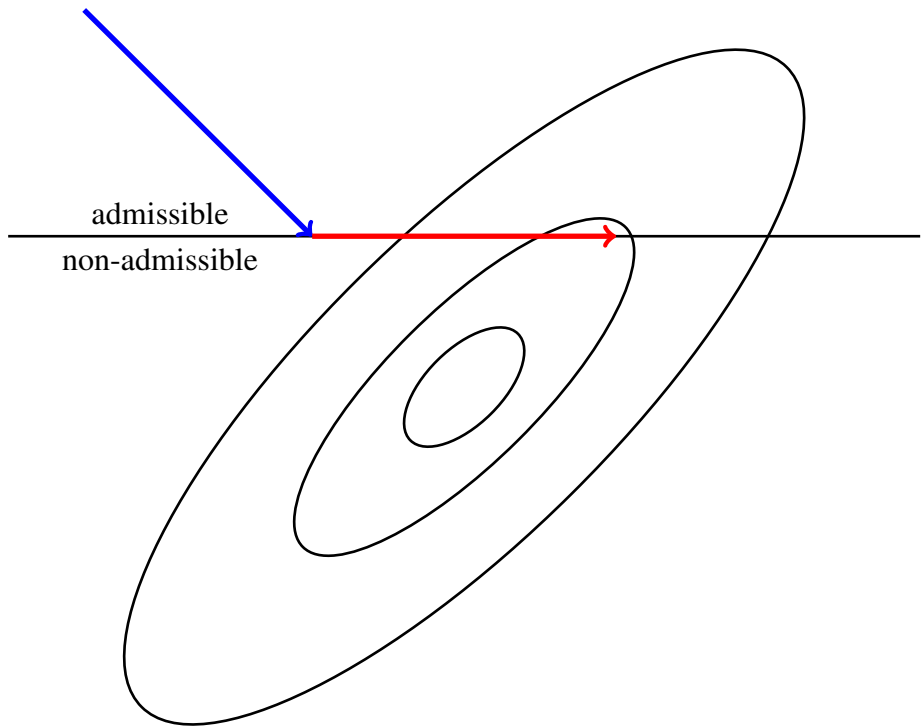
\begin{figure}
	\begin{tikzpicture}
		\node at (-4.0, 2.3) {admissible};
		\node at (-4.0, 1.7) {non-admissible};
		\draw [rotate=45, line width = 1pt] (0,0) ellipse (3cm and 1cm);
		\draw [rotate=45, line width = 1pt] (0,0) ellipse (6cm and 2cm);
		\draw [rotate=45, line width = 1pt] (0, 0) ellipse (1cm and 0.5cm);
		\draw [line width = 1pt] (-6.0, 2.0) -- (6.0, 2.0);
		\draw [blue, ->, line width = 2pt] (-5.0, 5.0) -- (-2.0, 2.0);
		\draw [red, ->, line width = 2pt] (-2.0, 2.0) -- (2.0, 2.0);
	\end{tikzpicture}
	\caption{Active set method in conjunction with the conjugate gradient method. The contour lines of the target function are depicted as ellipses, a line separates the admissible half-space from the non-admissible half-space. A blue and a red arrow are a first and second step of the active set method.  The blue step can be seen as a direct solution of the unrestricted problem using the CG method. As this step would leave the half-space of admissible solutions its length is restricted to the point where it collides with the hyperplane separating admissible and non-admissible solutions. The following step is constrained to the hyperplane, but no step length restrictions are needed. After two steps, the restricted minimizer is found.}
	\label{fig:ASMethod}
\end{figure}

The active set method is used to solve constrained minimization problems with inequalities as constraints \cite[section 16.4]{NW1999NO}
\[
	(c, d) = \argmin_{(c, d)} f(c, d) \quad \st \quad d_j \geq 0, \quad j = 1, \ldots, l.
\]
Our problem is obviously of this kind, as our convex set $C$, that forms our constraint, can be rewritten as a subspace of a vector space intersected with half-spaces in a Hilbert space. These half-spaces can be written in the form of inequalities, as above.
The method defines an active set, i.e. for every iterate $(c^k, d^k)$ a set of indices
\[
	\Aset_k\left(c^k, d^k\right) = \set{i}{d_i = 0}
\]
is defined that yields the inequality constraints that are strict in the point $\left(c^k, d^k\right)$. The task of the active set method is to move from iterate $\left(c^k, d^k \right)$ to the next iterate $\left(c^{k+1}, d^{k+1}\right)$ while respecting the set of constraints. This is done by limiting the step size to ensure that no change in the active set happens until the end of the step. This process is pictured in figure \ref{fig:ASMethod}.
In more formal terms the following steps are carried out:
\begin{enumerate}
	\item Start at a point $(c^1, d^1)$ that satisfies $ \forall j=1, \ldots, l:\quad d^1_j \geq 0.$
	\item Construct the active set 
		\[
			\Aset_k\left(c^k, d^k\right) = \set{i}{d^k_i = 0}.
		\]
	\item Solve
		\[
			(c, d) = \argmin_{(c, d)} \norm{\Aop (c, d) - u}_2 \quad \st \quad d_i = 0 \quad \forall i \in \Aset_k.
		\]
	\item Determine the largest $\lambda \in [0, 1]$ that satisfies
		\[
			\forall j = 1, \ldots, l: \quad d_j^{k + 1} \geq 0 \text{, with } \left(c^{k + 1}, d^{k + 1}\right) = (1-\lambda) \left(c^k, d^k\right) + \lambda (c, d).
		\] 
	\item Update the active set.
	\item Reiterate beginning with 3) until $\lambda = 1$ is satisfied in the previous step and no constraints can be dropped from the active set.
\end{enumerate}

\subsubsection{Determination of the Step Size}
	For the calculation of the step size the difference between a new solution of the quadratic program and the old solution is expressed as a direction
	\[
			\left(\Delta c, \Delta d \right) = \left(c^{k+1}, d^{k+1}\right) - \left(c^k, d^k\right).
	\]
	We search for the largest number $\lambda \in [0, 1]$ that satisfies
	\[
		\forall j=1, \ldots, l : \quad d^k_j + \lambda \Delta d_j \geq 0.
	\]
	For $\Delta d_j \neq 0$ we can express the $\lambda$ for which
	\[
		d^k_j + \lambda \Delta d_j = 0 \iff \lambda = - \frac{d^k_j}{\Delta d_j} 
	\]
	holds an set
	\[
		\lambda = \min\left(1,  \min_{j \in P} -\frac{d^k_j}{\Delta d_j}\right), \quad P = \set{j=1, \ldots, l}{\Delta d_j \neq 0 \wedge \frac{d^k_j}{\Delta d_j} > 0}.
	\]

\subsubsection{Updating the Active Set}
	A crucial step is the calculation of the active set
	\[
		\Aset_{k} \to \Aset_{k+1}
	\]
	 at the new position $\left(c^{k+1}, d^{k+1}\right)$  \cite[section 16.4]{NW1999NO}. 
	If the step size is smaller than $\lambda = 1$ a constraint, that was not in the active set, blocked the movement. This constraint must be part of the active set at the new position.
	If the active set allows no further movements it could be to large.  
	 An equality constraint can be deactivated if the steepest descent direction of the objective function points into the admissible region. At most one constraint should be deactivated during each update of the active set \cite{NW1999NO}.

\subsubsection{Solution of the Quadratic Program}
The restricted problem
\[
		(c, d) = \argmin_{(c, d)} \norm{\Aop (c, d) - \Aop u}_2 \quad \st \quad d_i = 0 \quad \forall i \in \Aset_k
\]
can be rewritten into a least squares problem 
\[
	(c, d) = \Nop \argmin_v \norm{\widetilde \Aop v - \Aop u}_2^2.
\]
Here, $\widetilde \Aop$ shall be an averaging operator with restricted basis 
\[
	\widetilde B = \sset{\phi_1, \ldots, \phi_k} \cup \set{\psi_i}{i \not \in \Aset_k}
\]
and the operator $\Nop$ shall embed zeros at the positions of the missing basis functions.
$\Aop$ and therefore $\widetilde \Aop$ have full column rank. Such least squares problems can be solved using the normal equations
\[
	w = \argmin_v \norm{\widetilde \Aop v - u}_2^2 \implies \underbrace{{\widetilde{\Aop}}^*\widetilde{\Aop}}_{O} w =  \widetilde{\Aop}^* u.
\]
The matrix $O$ of the normal equations is symmetric positive definite. The system can be solved via the Conjugate Gradiend method \cite{HS1952CG,Johnson2009Numerical}.
The following iteration
\[
	\begin{aligned}
		p_0 &= r_0 = k - A x_0, &
		a_i &= \frac{\skp{r_i}{r_i}}{\skp{p_i}{A p_i}} \\
		x_{i + 1} &= x_i + a_ip_i, &
		r_{i + 1} &= r_i - a_i A p_i \\
		b_i &= \frac{\skp{r_{i+1}}{r_{i+1}}} {\skp{r_i}{r_i}}, &
		p_{i+1} &= r_{i+1} + b_i p_i
	\end{aligned}
\]
is carried out. It is well known that this delivers the exact solution of a $n$ times $n$ equation system after $n$ iterations in exact arithmetic, and the convergence can be significantly faster for good initial estimates $x_0$. The conjugate gradient method was chosen as we solve similar linear systems during the active set method and can therefore reuse the last solution as initial value for the new iteration.

%% file: tests.tex
\subsection{Recovery Tests for Static Examples}

\begin{figure}
		\includegraphics[width=0.9\textwidth]{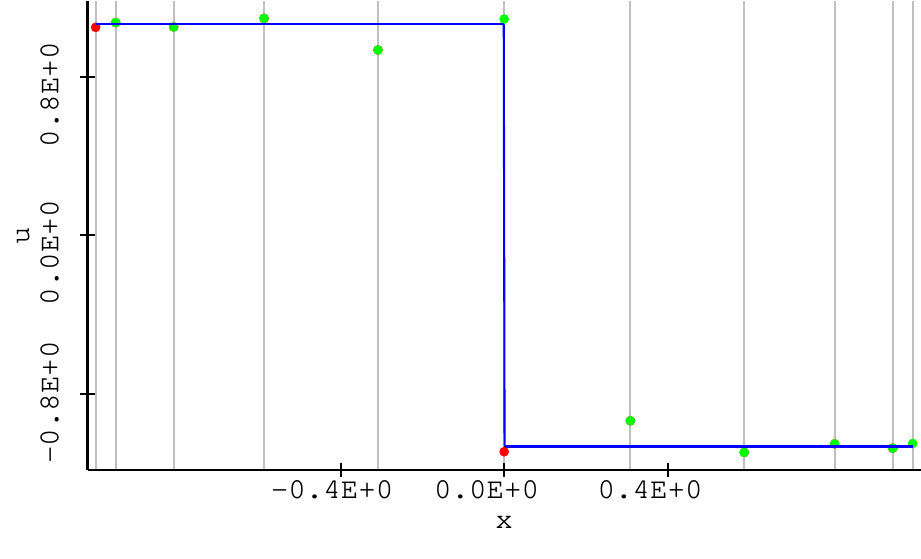}
		\caption{Initial condition $u_1$ with overlaid recovery. Cell edges are shown as vertical grey lines.
			The green circles are the right cell edge values, while the red circles give the left edge values of the recovered function. The jump height is recovered nearly exact. }
\label{fig:RecovTest1}
\end{figure}

\begin{figure}
		\includegraphics[width=0.9\textwidth]{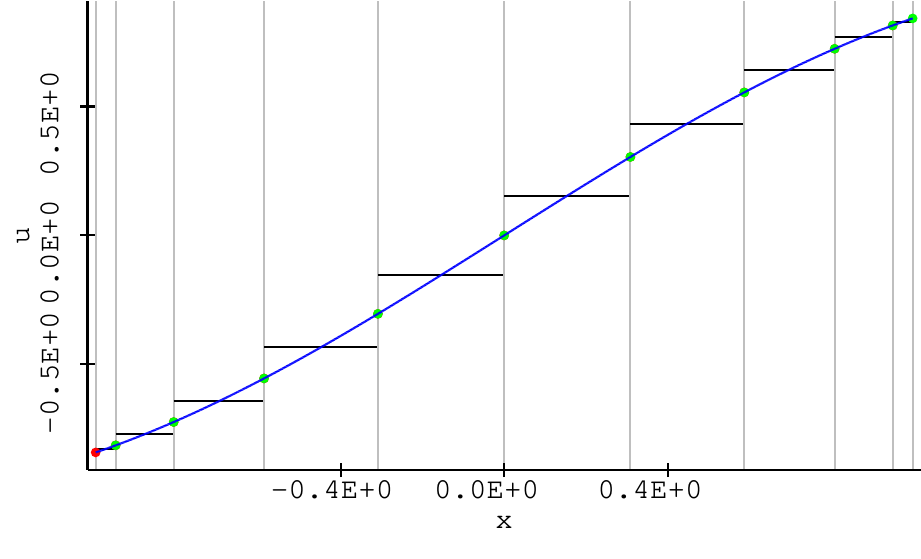}
		\caption{Recovery of a smooth function $u_2 = \sin(x)$ with overlaid recovery.  Cell edges are shown as vertical grey lines. The green circles are the right cell edge values, while the red circles give the left edge values of the recovered function. All function values are recovered with high accuracy. While two jumping basis functions are part of the admissible set of basis functions their contribution is negligible. The right edge value of every cell coincides with the left edge value of the neighbouring cell.}
		\label{fig:RecovTest2}
\end{figure}
\begin{figure}
		\includegraphics[width=0.9\textwidth]{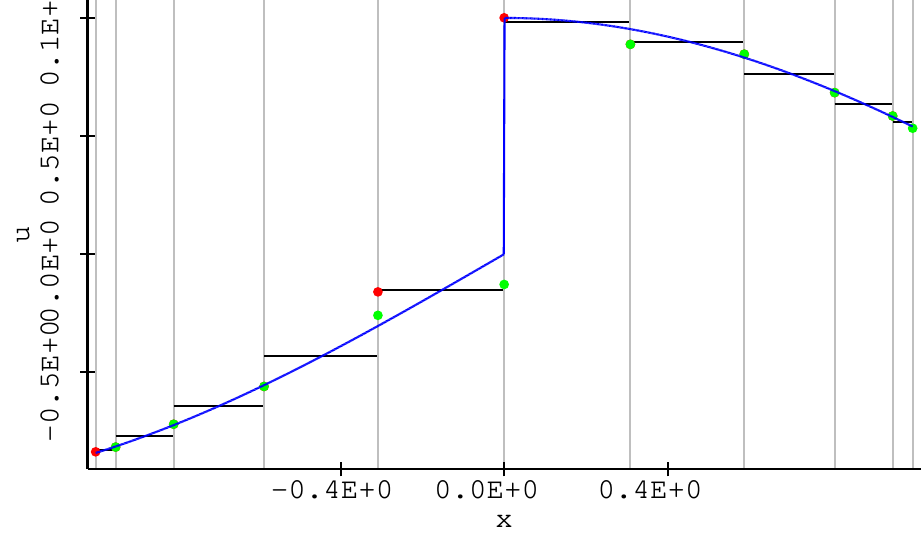}
		\caption{Recovery of function $u_3$ with smooth nonconstant areas connected by a discontinuity. Cell edges are shown as vertical grey lines. The green circles are the right cell edge values, while the red circles give the left edge values of the recovered function. As before, the jump height is predicted with high accuracy, the smooth parts are recovered with negligible oscillations.}
	\label{fig:RecovTest3}
\end{figure}

The recovery was constructed by enforcing simple constraints for the selected reconstructed function. These constraints are satisfied by classical ENO methods, but we did not enforce a ENO criterion, as there is no sensible definition of an essentially non-oscillative function. We will test if our assumption is true:
\nopagebreak \\
\emph{Adding jumps at the cell boundaries with the largest jumps in the averages and enforcing the sign criterion
allows for an Essentially Non-oscillatory reconstruction.}
\nopagebreak \\

We use the following functions:
\[
	u_1(x) = \begin{cases} 1 & x < 0 \\
							-1 & x \geq 0 \end{cases}, \quad
		u_2(x) = \sin(x), \quad 
		u_3(x) = \begin{cases} \sin(x) & x < 0 \\ 
								\cos(x) & x \geq 0 \end{cases}.
\]
The average values on a single macrocell from $-1$ to $1$ with $10$ subcells were calculated from these initial conditions, and we used our recovery algorithm with polynomials of degree less than or equal $8$ and $2$ variable jumps. The results are shown in Figures \ref{fig:RecovTest1}, \ref{fig:RecovTest2} and \ref{fig:RecovTest3}.

The first function is recovered with nearly perfect resolution of the jump height. Slight oscillations are visible that vanish with growing distance from the discontinuity at $x = 0$.
Our second example shows that a smooth function is recovered without any jumps. 
The third example combines a jump with smooth but non-constant areas. We note that because our recovery algorithm is non-linear, the result of such a recovery is not clear from the first two examples.
In this case two discontinuities are recovered while the true function has only one discontinuity. Still, the second recovered discontinuity respects the sign property and oscillations are negligible.

\subsection{Tests of the Recovery in Finite-Volume Methods }
Now, we will test our new recovery in a FV method for the Euler equations of gas dynamics. 
The flux function in equation \eqref{eq:HCL} is given by \cite{Harten1983On}
\[
	f(u) = \begin{pmatrix}
	\rho v \\ 
	\frac{\rho  v^2}{2} + p \\
	v (E + p) 
	\end{pmatrix}, \quad u = \begin{pmatrix}
	\rho \\ 
	\rho v \\
	E
	\end{pmatrix}, \quad p = (\gamma -1) \left( E  - \frac {\rho v^2}{2} \right).
\]
 The HLL flux
 \cite{HLL1983On} 
 \[
 	f(u_l, u_r) = \begin{cases} f(u_l) & 0 < a_l \\
 								\frac{a_r f(u_l) - a_l f(u_r) + a_r a_l (u_r-u_l)}{a_r - a_l} & a_l < 0 < a_r \\
 								f(ur) & a_r < 0
 								\end{cases} , \quad 
 \]
 will be used with the simple speeds from Davis \cite{Davis1988Simplified}
 \[
 	a_l = \min(v_l - c_l, v_r - c_r), \quad a_r = \max(v_l + c_l, v_r + c_r).
 \]
 Here, $v_l$ and $v_r$ are the fluid velocities on the left and right of the discontinuity, while $c_l$ and $c_r$ are the sound speeds on the left and right of the discontinuity.
 
 Our method has three integer parameters: The number of subcells per macrocell $q+1$, the number of continuous basis functions $k$ and the number of jumps $l$. The basic compatability relation $k + l \leq q+1$ restricts the possible choices of these parameters. We propose to choose $q+1 = 4$ or $q+1 = 8$ as this allows to fit the resulting matrices and vectors into SIMD registers. The number of jumps was set to $l = 1$. To achieve the highest possible accuracy for smooth solutions $k = q+1 - l$ smooth basis functions were used. Experiments were carried out for $q+1 = 2, \ldots, 11$ and the method worked as expected. We will only report the results for $q+1  = 4$ and $q+1=8$. 
 
 Time integration will be carried out using the SSPRK(3, 3) method with the time-step controlled via a CFL-type restriction
 \[
 	\Delta t = c_{fl} \frac{\Delta x}{c_{\mathrm{max}}}.
 \]
 The constant $c_{fl}$ is set to $ c_{fl} = 1/10$, the spacial grid size $\Delta x$ is the minimal size of a subcell in the grid and $c_\mathrm{max}$ is the highest signal speed of the solution.
 \begin{enumerate}
 	\item 
 Our first example is a numerical convergence experiment using the initial condition
 \[
 	\rho(x, 0) = 1 + \e^{-(x-1)^2/2}, \quad v(x, 0) = 1, \quad p(x, 0) = 1
 \]
with periodic boundary conditions. This results in an advection of a density variation from the left to the right with the analytic solution
\[
	\rho(x, t) = 1 + \e^{-(x-t-1)^2/2}, \quad v(x, t) = 1, \quad p(x, t) = 1.
\]
We use a periodic grid with $N \in [16, 20, 24, 28, 32, 36, 40, 44, 48, 52]$ macrocells and compare the numerical solution to the exact solution.
The resulting errors are plotted in Figure \ref{fig:NumConv}. If $k$ monomials are used as the smooth ansatz functions one expects a $k$-th order convergence, as the recovery polynomials are of degree $k-1$. The numerical tests indicate a $\txfrac  5 2$ order convergence for a scheme using $k = 3$ continuous basis functions. A scheme using $k=7$ basis functions converges with order six. The reason for this degradation is unclear, but one can conjecture that a convergence of order 7 can be observed at extreme resolutions that are not of practical interest. Further, one can conjecture that different subcell layouts could restore the convergence speed at lower resolutions. 

\begin{figure}
	\begin{subfigure}{0.99 \textwidth}
	\includegraphics[width=0.9\textwidth]{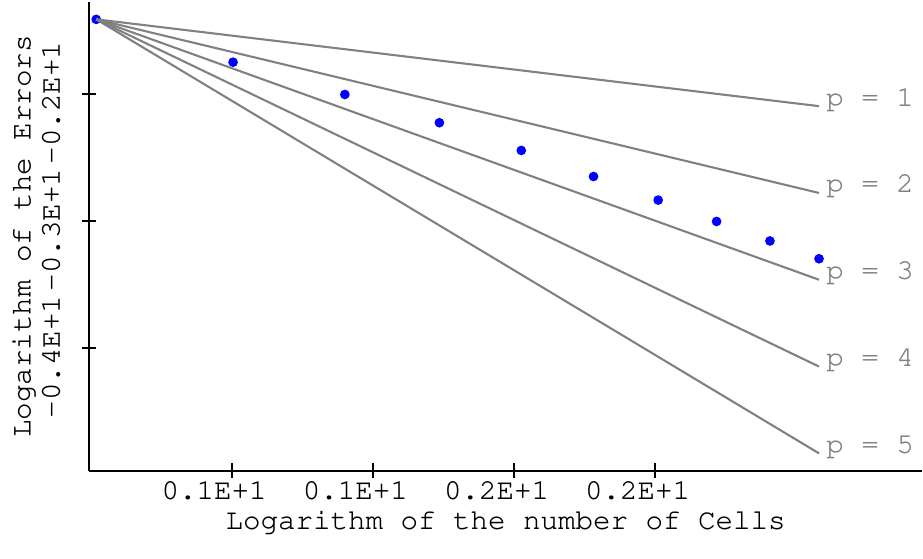}
		\caption{Errors of a numerical solution calculated using $4$ subcells, $3$ continuous basis functions and one jump per subcell.}
	\end{subfigure}
		\begin{subfigure}{0.99 \textwidth}
		\includegraphics[width=0.9\textwidth]{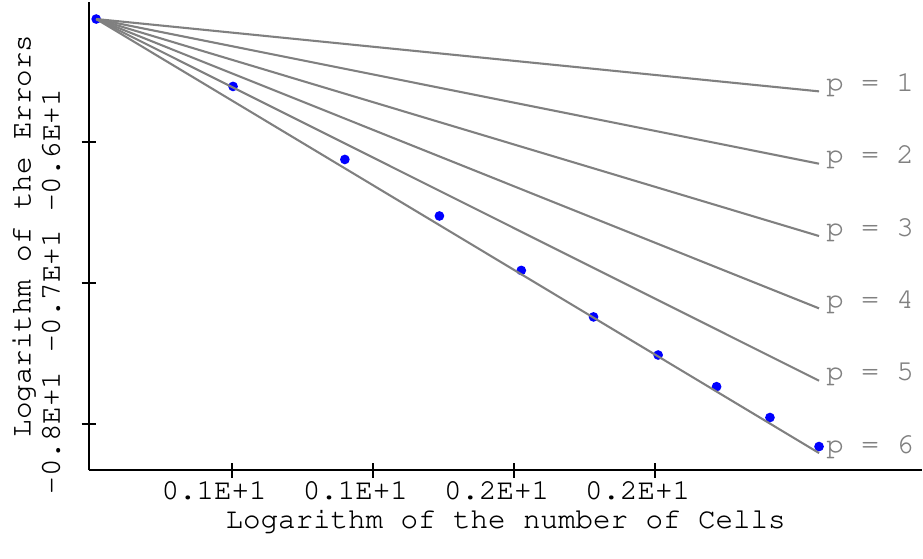}
		\caption{Errors of a numerical solution calculated using $8$ subcells, $7$ continuous basis functions and one jump per subcell.}
	\end{subfigure}
	\caption{Numerical convergence analysis using a density variation advection up to $t=10$.}
	\label{fig:NumConv}
\end{figure}

\item As second initial condition Sods shock tube \cite{Sod1978Survey} 
\[
	u(x, 0) = \begin{cases} \rho = 1.000, \quad v= 0, \quad p = 1.0 & x \leq 0, \\ \rho = 0.125, \quad v = 0, \quad p = 0.1 & x \geq 0. \end{cases}
\]
will be used.
\begin{figure}

		\includegraphics[width=\epswidth]{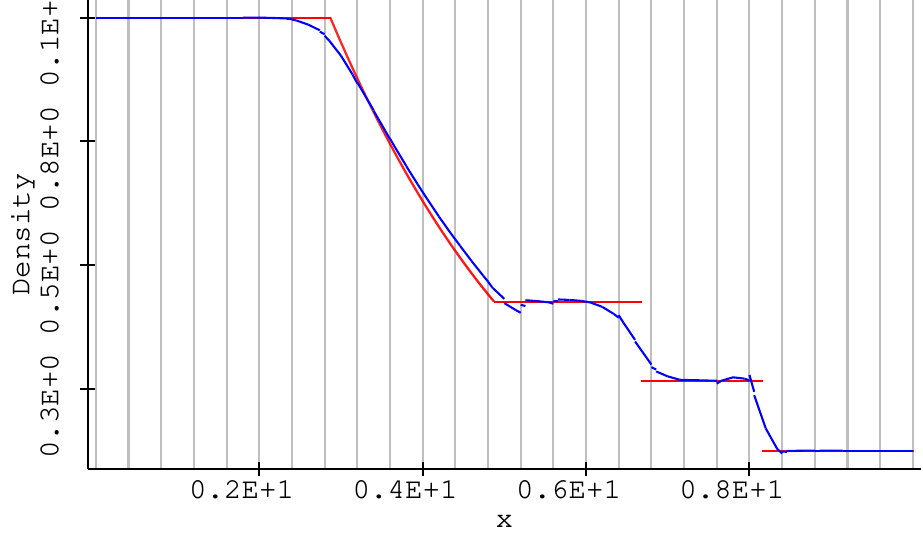}
	\caption{Solution to Sods' problem. A total of $N = 100$ cells were used. These were grouped into $K = 25$ cell groups of $4$ cells.
		Quadratic polynomials and one discontinuity per cell group were used in the recovery procedure. The density is shown at $t = 1.8$.}
	\label{fig:sod4}
\end{figure}
\begin{figure}
		\includegraphics[width=\epswidth]{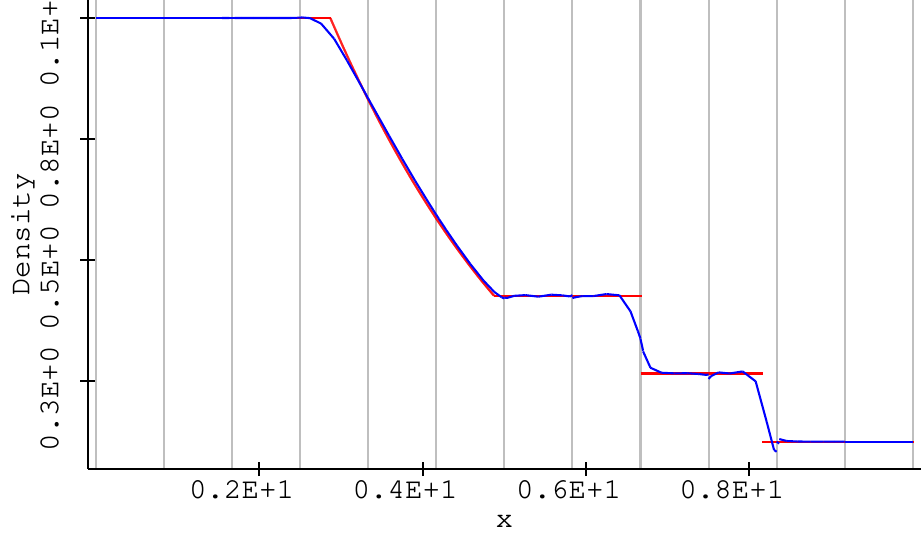}
	\caption{Solution to Sods' problem. A total of $N = 96$ cells were used. These were grouped into $K = 12$ cell groups of $8$ cells.
		Polynomials of degree $6$ and one discontinuity per cell group were used in the recovery procedure. The density is shown at $t = 1.8$.}
	\label{fig:sod8}
\end{figure}
The results can be seen in figure \ref{fig:sod4} and \ref{fig:sod8}. The Riemann problem to this initial condition consists of a left rarefaction, a right moving contact disconinuity and a faster right moving shock. Our reference solution is given by an exact Riemann solver following Toro's approach \cite{Toro2009Riemann}. Let us remark that both discontinuities lie in the middle of macrocells in figure \ref{fig:sod4}, but both discontinuities are represented via continuous high gradient polynomials. The shock lies on a macrocell boundary in figure \ref{fig:sod8} but is still represented without a discontinuity. Albeit the solution is non-smooth the errors of the solution calculated with a higher order method are lower.

\item The Riemann problem with the initial condition 
\[
	u(x, 0) = \begin{cases} \rho = 0.445,\quad v = 0.698,\quad p = 3.528 & x < 0, \\ \rho = 0.500,\quad v = 0.000,\quad p = 0.571& x > 0.\end{cases}
\]
is a common test \cite{SO1989,HEOC1987UniformlyIII}. We will refer to it as Lax's problem.
\begin{figure}
		\includegraphics[width=\epswidth]{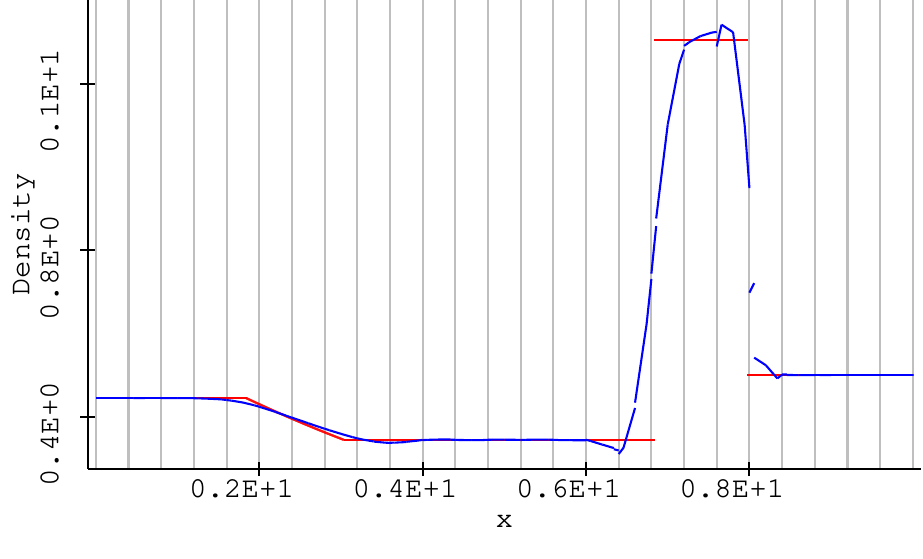}
	\caption{Solution to Lax problem. A total of $N = 100$ cells were used. These were grouped into $K = 25$ cell groups of $4$ cells.
		Quadratic polynomials and one discontinuity per cell group were used in the recovery procedure. The density is shown at $t = 1.2$.}
	\label{fig:lax4}
\end{figure}
\begin{figure}
		\includegraphics[width=\epswidth]{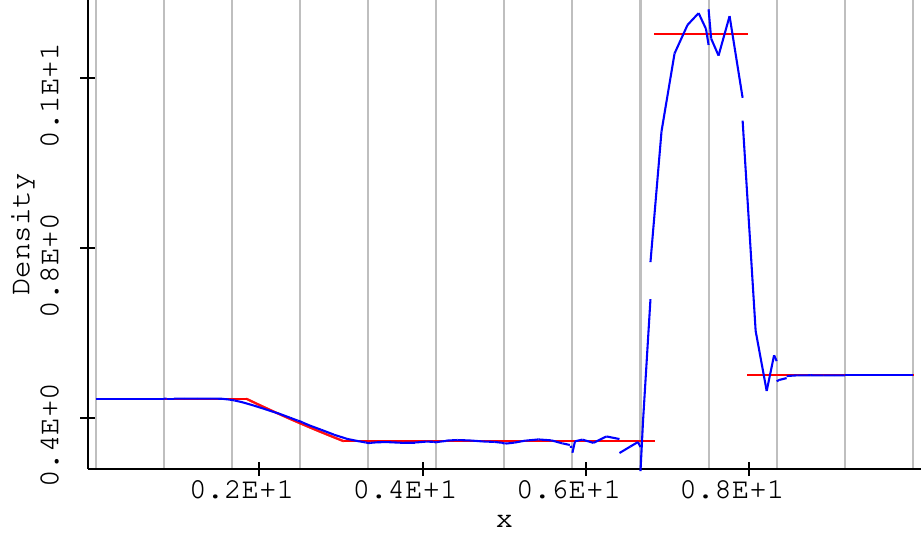}
	\caption{Solution to Lax problem. A total of $N = 100$ cells were used. These were grouped into $K = 25$ cell groups of $4$ cells.
		Polynomials of degree $6$ and one discontinuity per cell group were used in the recovery procedure. The density is shown at $t = 1.2$.}
	\label{fig:lax8}
\end{figure}
The results at $t = 1.2$ can be seen in figures \ref{fig:lax4} and \ref{fig:lax8}. 
Once more the solution consists of a left moving rarefaction, and a right moving contact and shock.
Here, discontinuities in macrocells are represented by discontinuous functions, but the jumps in these functions are significantly smaller than the true jumps. 
We can only assume that this is due to the fact that the functions are recovered from approximate average values that were calculated by a numerical scheme. The inherent dissipation of the numerical scheme smoothed out the discontinuities until only smaller discontinuities or high gradients were left.

\item A classical testcase for the ability of a recovery is the Shu-Osher testcase given by the initial condition
\[
	u(x, 0) = \Bigg \{ \begin{array}{l l l l}
		\rho = 3.857143,& v = 2.629369, & p = 10.33333& x < 1, \\
		\rho = 1 + \epsilon \sin(5 x), & v = 0.0,&  p = 1.0& x > 1.
	\end{array}
\]
with $\epsilon = 0.2$ \cite{SO1989}.
The results can be seen in figures \ref{fig:so4LR}, \ref{fig:so4HR}, \ref{fig:so8LR} and \ref{fig:so8HR}.
\begin{figure}
		\includegraphics[width=\epswidth]{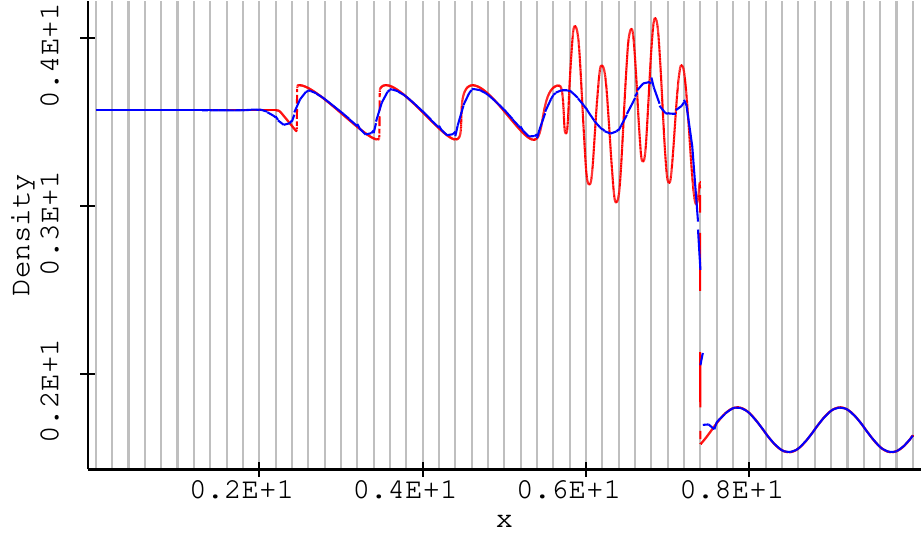}
	\caption{Solution to the Shu-Osher test. A total of $N = 200$ cells were used. These were grouped into $K = 50$ cell groups of $4$ cells.
		Quadratic polynomials and one discontinuity per cell group were used in the recovery procedure. The density is shown at $t = 1.8$.}
	\label{fig:so4LR}
\end{figure}
\begin{figure}
		\includegraphics[width=\epswidth]{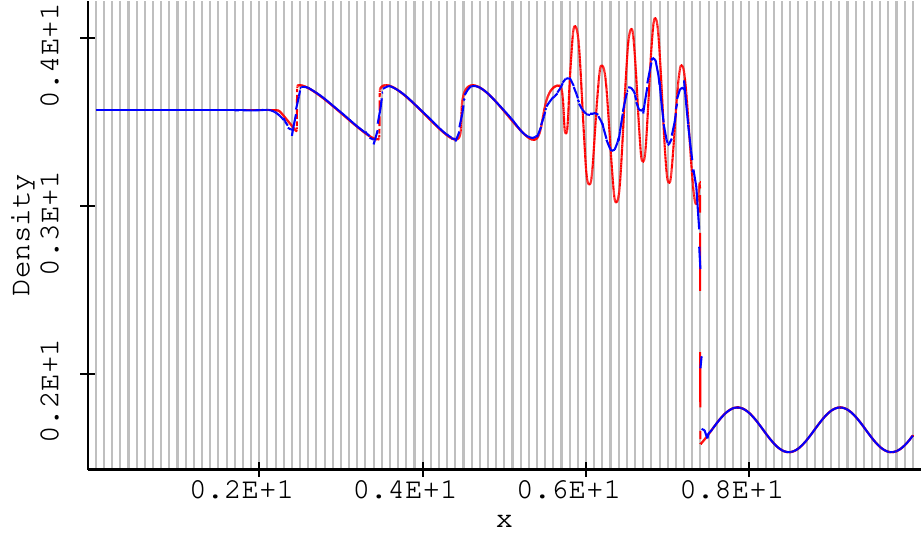}
\caption{Solution to the Shu-Osher test. A total of $N = 400$ cells were used. These were grouped into $K = 100$ cell groups of $4$ cells.
Quadratic polynomials and one discontinuity per cell group were used in the recovery procedure. The density is shown at $t = 1.8$.}
	\label{fig:so4HR}
\end{figure}
\begin{figure}
		\includegraphics[width=\epswidth]{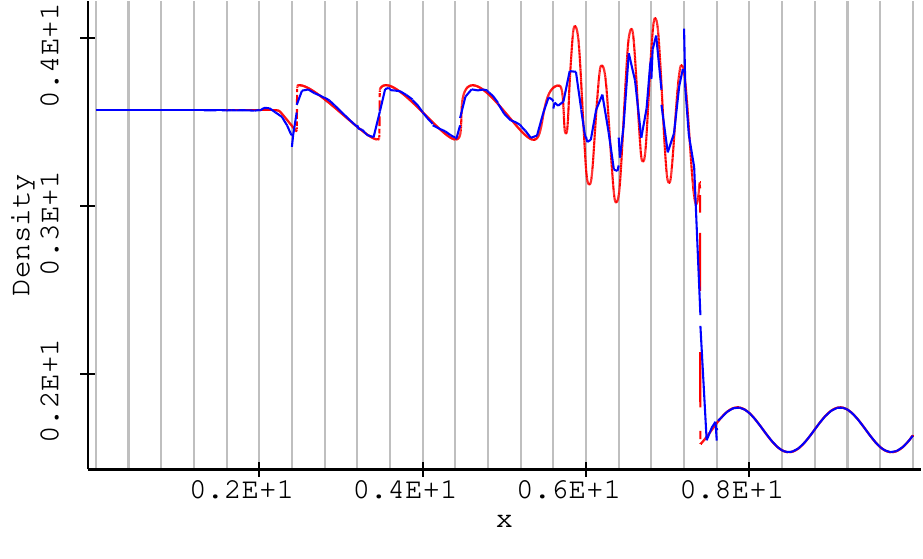}
	\caption{Solution to the Shu-Osher test. A total of $N = 200$ cells were used. These were grouped into $K = 25$ cell groups of $8$ cells.
		Polynomials of degree $6$ and one discontinuity per cell group were used in the recovery procedure. The density is shown at $t = 1.8$.}
	\label{fig:so8LR}
\end{figure}
\begin{figure}
		\includegraphics[width=\epswidth]{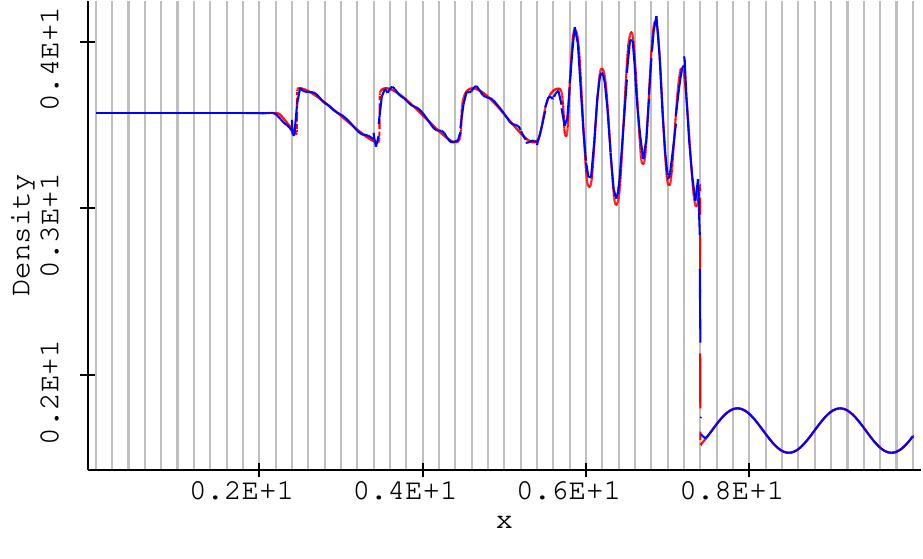}
	\caption{Solution to the Shu-Osher test. A total of $N = 400$ cells were used. These were grouped into $K = 50$ cell groups of $8$ cells.
		Polynomials of degree $6$ and one discontinuity per cell group were used in the recovery procedure. The density is shown at $t = 1.8$.}
	\label{fig:so8HR}
\end{figure}
The reference solution was calculated by a MUSCL scheme with $N = 8192$ cells and the HLL numerical flux. If a scheme with $4$ subcells per macrocell and $50$ macrocells is used a total of $200$ cells correspond to the classical number of cells for this test-case. The resulting solution in figure \ref{fig:so4LR} is only a rough approximation of the true solution and the low wavelength waves around $x = 6.5$ are missing. If $400$ cells are used, as in figure \ref{fig:so4HR}, the approximate solution predicts some short wavelength waves.
The scheme using $8$ subcells has a significantly higher resolution using $200$ cells, consult figure \ref{fig:so8LR}. The fine structure around $x = 6.5$ is visible with all maxima and minima.
The approximate solution using $400$ cells is indistinguishable from the exact solution.

\end{enumerate}

%% file: concl.tex
We developed a non-linear reconstruction procedure for Spectral Volume (SV) methods. 
The recovered function satisfies the sign property \cite{FMT2013ENO} for all jumps inside spectral volumes.
A key ingredient was the selection of suitable basis functions that can represent jumps that are not located at the edges of Spectral Volumes.
The system of linear equations employed in SV methods was generalized to a least squares system of equations with additional convex constraints.
The active set method was employed to solve these optimization problems efficiently.
Numerical tests indicate that the recovered functions oscillate in the pressence of discontinuities, but their oscillations are small compared to the jump heights. 
Further tests in Spectral Volume solvers indicate high order accuracy and essentially oscillation free solutions as the oscillations are confined to the shocked macrocells and the heights of the oscillations in the shocked cells are small compared to the jump heights.
A sequel to this publication will show how the recovery technique can be used on two-dimensional grids.